\documentclass[review,11.5pt]{elsarticle}
\usepackage{ifpdf}
\usepackage{graphicx,natbib,amssymb,lineno, latexsym, theorem,amsmath}

\ifpdf
\usepackage[%
  pdftitle={Dynamics of the Non-autonomous Owen-Smith Model},%
  pdfauthor={Mohamed A. A. Bakheet},%
  pdfsubject={},%
  pdfkeywords={metaphysiological, permanence, extinction, asymptotic stability, periodic solutions},%
  pdfstartview=FitH,%
  bookmarks=true,%
  bookmarksopen=true,%
  breaklinks=true,%
  colorlinks=true,%
  linkcolor=blue,anchorcolor=blue,%
  citecolor=blue,filecolor=blue,%
  menucolor=blue,pagecolor=blue,%
  urlcolor=blue]{hyperref}
\else
\usepackage[%
  breaklinks=true,%
  colorlinks=true,%
  linkcolor=blue,anchorcolor=blue,%
  citecolor=blue,filecolor=blue,%
  menucolor=blue,pagecolor=blue,%
  urlcolor=blue]{hyperref}
\fi

\makeatletter
\def\elsartstyle{%
    \def\normalsize{\@setfontsize\normalsize\@xiipt{14.5}}
    \def\small{\@setfontsize\small\@xipt{13.6}}
    \let\footnotesize=\small
    \def\large{\@setfontsize\large\@xivpt{18}}
    \def\Large{\@setfontsize\Large\@xviipt{22}}
    \skip\@mpfootins = 18\p@ \@plus 2\p@
    \normalsize
}
\@ifundefined{square}{}{}
\makeatother

\newcommand{\RR}{\mathrm{I\!R\!}}
\theoremheaderfont{\upshape} 
\theoremheaderfont{\itshape} {\theoremstyle{break}

\newtheorem{Lem}{\textbf{Lemma}}[section] \theoremstyle{break}
\newtheorem{Thm}{\textbf{Theorem}}[section] {\theoremstyle{break}
\theorembodyfont{\rmfamily} 
{\newtheorem{Prf}{\textbf{Proof}}[section]}
\newtheorem{Def}{\textbf {Definition}}[section]
{\theoremstyle{plain}
  \theorembodyfont{\rmfamily} 
}}   
 
\begin{document}

\begin{frontmatter}
\title{Dynamics of the Non-autonomous Owen-Smith Model}
\author{Mohamed A. A. Bakheet}
\ead{ mohamed.bakheet@uct.ac.za}
\author{Henri Laurie}
\address{Department of Mathematics and Applied Mathematics,\\ University of Cape Town, Rondebousch 7701\\ South Africa}

\begin{abstract}
In this paper we study the dynamics of the general case of Owen-Smith metaphysiological model, to explore the effects of seasonality on population fluctuations. The study will include the permanence, herbivore extinction, global asymptotic stability and existence of positive periodic solutions. Under certain assumptions, we obtained sufficient and necessary conditions which guarantee the permanence of herbivore and vegetation species, and existence of periodic solutions. The techniques we used here are, comparison method of differential equations, Lyapunov method and Brouwer's fixed-point Theorem.
\end{abstract}

\end{frontmatter}

\section{Introduction}
Owen-Smith in \cite{Smith04} presents the following metaphysiological model (a metaphysiological relates aggregated population dynamics to biomass gained from resources consumed, relative to physiological attrition and mortality losses \cite{Smith2002a,Getz1991,Getz1993}), to explore how resource heterogeneity functions to dampen herbivore population fluctuations
\begin{eqnarray}\label{smith}
\left. 
\begin{array}{l}
\frac{dV}{dt} = r V\left(1 \,-\, \frac{V}{K}\right)- \frac {i_m(V \,-\, v_u)} {b_i \,+\, V \,-\, v_u} H \\
\frac{dH}{dt} = H\left( \frac{Ci_m(V \,-\, v_u) }{b_g \,+\, V \,-\, v_u} - m_p - q_0 - q_s -\frac{qm_p (b_g \,+\, V \,-\, v_u)}{Ci_m(V \,-\, v_u)} \right)
\end{array} \right.
\end{eqnarray}
where $V,H$ represent vegetation, herbivore densities respectively, $r$ is the maximum relative growth rate of vegetation at time $t$, $K$ is the vegetation carrying capacity, $i_m$ is the maximum intake rate by unit herbivore, $b_i$ is the vegetation biomass at which herbivore intake rate reaches half of its maximum, $b_g$ is the half-saturation level which is the vegetation biomass at which the intake rate becomes half of its maximum, $v_u$ is the ungrazable amount of vegetation at which the intake rate of the herbivores becomes zero (vegetation reserve), $C$ is the rate of  conversion from consumed vegetation biomass into herbivore biomass, $m_p$ is the relative rate of physiological attrition in herbivore biomass, $q_s$ is the basic rate of loss of herbivore biomass through mortality due to senescence, and $q_0$ determine how the steepness of nutritionally-related mortality increases with diminishing food gains. For more background on system \eqref{smith} see \cite{Smith04,Smith2002a,Smith98,Smith02b}.
\\ 

Most environmental factors are highly variable and they show seasonal variations, and in response, birth rates, death rates and other vital rates of populations vary greatly in time \cite{Fan}. In order to represent the effects of seasonality in his model, Owen-Smith \cite{Smith04} made the vegetation growth to be seasonal by using a step production function assuming that vegetation growth occurred continually during growing season, and then ceased during dormant season (i.e. $r = r_v$ for the growing season weeks, and $r = 0$ for the dormant season weeks). Seasonality enter Owen-Smith model \eqref{smith} only via the maximum vegetation growth rate $r$. However food abundance and atmospheric temperature show high variation in between seasons, animal intake rates are responses to these factors. The vegetation reserve represents the amount of vegetation that herbivores cannot eat, and this amount may also show distinct variations during wet and dry seasons. Herbivore death rate is caused by many environmental factors such as predation, food abundance and food quality, and all these factors respond to seasonal changes. So, it is plausible to assume that all the parameters in the Owen-Smith model \eqref{smith} change with time.
\\

To perform our analysis firstly, assume that all the model parameters are continuous $\omega$-periodic functions. Secondly, assume the half saturation rates for consumption and conversion are equal, i.e $ b_i\,=\, b_g$. Thirdly, and for biological realism, we assume that $v(t_0) \,>\, 0 \,\, \text{and} \,\, h(t_0) \,>\, 0$. Finally, we simplify system \eqref{smith} by making change of variables. Let $H=h,\, V=v, \, a = r_v,\, b=\frac{r_v}{K},\, c=i_m,\, \alpha=Ci_m,\, \beta_1= b_i- v_u,\, \beta_2 =  b_g - v_u,\, \gamma=\frac{qm_p}{Ci_m}, \, \rho \,=\, v_u$, and $R= m_p+q_0+q_s$. With these assumptions the model given by \eqref{smith} becomes the following:
\begin{eqnarray}\label{NAM}
\left. 
\begin{array}{l}\label{smith2}
\frac{dv}{dt} = v(t) \left( a(t) - b(t)v(t) \right)- c(t) \frac { v(t) \,-\, \rho(t)} {\beta(t) \,+\, v(t)}h(t)
\\
\frac{dh}{dt} = h(t) \left( \alpha(t) \frac{v(t) - \rho(t) }{\beta(t) \,+\, v(t)} - R(t) - \gamma(t) \frac{\beta(t) \,+\, v(t)}{v(t) \,-\, \rho(t)} \right)
\end{array} \right.
\end{eqnarray}
where $a(t),b(t),c(t),\beta(t),\alpha(t), \gamma(t), \rho(t)$, and $R(t)$ are all continuous and positive $\omega-$periodic functions. \\ We developed the research method given in \cite{Fan,Yang2008,Chen2006} and \cite{Wang}, and we established new sufficient conditions for which herbivore and vegetation having globally asymptotically solutions and having periodic solutions. Also we established new necessary condition for which species for system \eqref{NAM} remain permanent.
 
\section{Preliminaries}
Before we going to our main results, we give some definitions, notations, lemmas and theorems that we are going to use in order to obtain the main results.
\begin{Def}
For a positive continuous and $\omega$-periodic function $f(t)$, we denoted by $A_\omega (f)$ the average of $f(t)$ over an interval of length $\omega$
\begin{equation}\label{average}
A_\omega(f) = \frac{1}{\omega×} \int_{0}^\omega f(t) dt. 
\end{equation}
\end{Def}
\begin{Def}
For a bounded continuous function $g(t)$ on $\RR$ we use the following notations:
\begin{equation}
g^u \,=\, \sup_{t\in \RR} \, g(t), \quad  g^l \,=\, \inf_{t\in \RR} \, g(t). \nonumber
\end{equation}
\end{Def}
The following theorem proves the positive invariant of any solutions to \eqref{NAM} with positive initial vales. Let $\RR_{+} \,=\, [ 0, \infty) $,  and $\RR_{+}^{N} \,=\, \left\{ X \in \RR^{N} | \, X \, \geq \, 0  \right\}$.
\begin{Thm}
If $\left( v, h \right) \in \RR_{+}^{2}$ is a solution to \eqref{NAM} with positive initial values, then both $\left( v(t), h(t) \right)$ will remains positive for all $t \geq t_0$.
\end{Thm}
\begin{Prf} 
By integrating both side of system \eqref{NAM} with respect to $t$ from $t_0$ we obtain the following:
\begin{eqnarray}
v(t) &=& v(t_0) \exp\left[\int_{t_0}^t \left(  a(s)-b(s)v(s)  - \frac {c(s) \left( v(s) - \rho (s) \right) h(s)} {v(s)\left( \beta(s)+v(s) \right)} \right)ds \right] \nonumber \\
h(t) &=& h(t_0) \exp\left[\int_{t_0}^t \left( \frac{\alpha(s)( v(s) - \rho (s) )}{\beta(s)+v(s)} - R(s) -\gamma(s)\frac{\beta(s)+v(s)}{v(s) - \rho(s)} \right) ds \right] \nonumber
\end{eqnarray}
It is obvious when $v(t_0)>0$, and $h(t_0)>0$, both $v(t)$ and $h(t)$ are positive and then invariant for all $t \geq t_0$. The proof is complete.
\end{Prf}
\begin{Def}[\cite{Wang}]\label{ultimately-bounded-def} 
The solution of system \eqref{NAM} is said to be ultimately bounded if there exists $B>0$ such that for every solution $\left(v(t),h(t)\right)$ of system \eqref{NAM} there exist $T>0$ such that $\| \left(v(t),h(t)\right)\| \leq B$, for all $ t\geq t_0+T$, where $B$ is independent of particular solutions while $T$ can depend on the solution.
\end{Def}
\begin{Def}[\cite{Liu2003}]
The system \eqref{NAM} is said to be permanent if there exists a compact set $\Gamma$ in the interior of $ \RR_+^{2} = \left\{ (v,h) \in \RR^2 | v \geq 0, h \geq 0 \right\} $ such that all the solutions starting in the interior of $\RR_+^{2}$ ultimately enter $\Gamma$ and remains in it. In other words the system is permanent if there exists positive constants $\delta$, and $\Delta$ with $0< \delta <\Delta$ such that
\begin{align}\label{per}
 & \min \left[ \lim_{t\rightarrow +\infty} \inf v(t), \lim_{t\rightarrow +\infty} \inf h(t) \right]  \geq \delta \nonumber \\
 & \max \left[ \lim_{t\rightarrow +\infty} \sup v(t), \lim_{t\rightarrow +\infty} \sup h(t) \right]  \leq \Delta \nonumber 
\end{align}
for all solutions of system \eqref{NAM} with positive initial values.
\end{Def}
\begin{Def}[\cite{Meng,Wang,Xia2008}] \label{globally-stable}
A bounded non-negative solution $(\hat v (t), \hat h(t))$ of \eqref{NAM} is said to be globally asymptotically stable (or globally attractive) if, for any other solution $(v(t), h(t))$ of \eqref{NAM} with positive initial values, the following condition holds:
\begin{equation}\label{gstable}
\lim_{t\rightarrow +\infty} \left( | v(t) - \hat v (t) | + | h(t) - \hat h (t) | \right) \,=\, 0. \nonumber
\end{equation}
\end{Def}
The following Lemma proves that, system \eqref{NAM} has a unique positive globally asymptotically $\omega$-periodic solution $\left( v^{*},0 \right)$. Using this theorem we can prove positive invariance, permanence, and global asymptotic stability of the full model \eqref{NAM}.
\begin{Lem}\label{globally-stable-lemma}
 If $a(t)$ and $b(t)$ are $\omega$-periodic, and if $A_\omega(b)>0$ and $A_\omega(a) >0$, then system \eqref{NAM} has a unique positive $\omega$-periodic solution $\left( v^{*},0 \right)$ which is globally asymptotically stable, where $A_\omega(\cdot)$ defined by \eqref{average}, and 
\begin{eqnarray}\label{global-asymptotic-solution}
 v^{*}(t)= \exp \left\{\int_{0}^t a(s)ds \right\} \left[ \int_{0}^t b(s)\exp \left\{ \int_{0}^s a(\nu) d\nu \right\} ds + c^*\right]^{-1}.
\end{eqnarray}
\end{Lem}
\begin{Prf}
When $h=0$, system \eqref{NAM} is reduced to the following differential equation:
\begin{equation}\label{gas}
 \dot v(t) =v(t) \left[a(t)-b(t)v(t) \right].
\end{equation}
It is not difficult to show that \eqref{global-asymptotic-solution} is $\omega$-periodic solution to \eqref{gas}. To prove $v^*(t)$ is asymptotically stable consider the Lyapunov function
\begin{eqnarray}\label{W}
W(t) \,=\, \left| \ln\{v(t)\}-\ln\{v^*(t)\} \right|
\end{eqnarray}
where $v(t)$ is any positive bounded solution to \eqref{gas}. Calculating the upper right derivative of $W(t)$ along the solution of \eqref{gas}, we obtain:
\begin{eqnarray}\label{lyap1}
\dot W (t) &=& \text{sgn} \left\{v(t) - v^*(t) \right\}\left[ \frac{\dot v(t)}{v(t)} - \frac{\dot {(v^*)}(t) }{v^*(t)} \right]
\end{eqnarray}
and then we have
\begin{eqnarray}\label{lyap11}
\dot W(t) &=& \text{sgn}\left\{v(t) - v^*(t) \right\} \,\, \left[ a(t)-b(t)v(t)- a(t) + b(t)v^*(t)) \right] \nonumber\\ 
&=& \text{sgn} \left\{v(t) - v^*(t) \right\} \,\, \left[-b(t)v(t) + b(t)v^*(t) \right] \nonumber\\
&=& -b(t)\,\, \text{sgn} \left\{v(t) - v^*(t) \right\} \,\, \left[v(t) - v^*(t) \right] \nonumber\\
&=&-b(t)\,\,|v(t)-v^*(t)|
\end{eqnarray}
Assume that $\zeta= \max \left\{ |\ln\{v(t)\}|, |\ln \{v^*(t)\}| \right\} $.
Let $ \zeta^* = \exp\left\{ \zeta \right\}$, then we have:
 \begin{eqnarray}
|v(t)-v^*(t)| \leq \zeta^* \, W(t) \nonumber
\end{eqnarray}
substituting this back into \eqref{lyap11} we obtain:
\begin{eqnarray}\label{lyap2}
\dot W(t) \leq - b(t) \zeta^* \, W(t)
\end{eqnarray}
integrating both side of \eqref{lyap2} from $t_0$ to $t$ yields:
\begin{eqnarray}
W(t) \leq W(t_0) \exp\left\{ -\zeta^* \int_{t_0}^t b(s)ds \right\}.
\end{eqnarray}
And since $A_\omega\left( b \right) \,>\, 0$, then $ \lim_{t \rightarrow \infty} W(t)=0$. Thereby $ \lim_{t \rightarrow \infty} |v(t)-v^*(t)|=0$.
The proof is complete.   
\end{Prf}
\begin{Thm}\label{M_1-M_2-theorem}
If $\left( v(t), h(t) \right)$ is any positive solution to \eqref{NAM}.
Then there exists positive constants $M_1$ and $M_2$ such that
\begin{eqnarray}
 \lim_{t \rightarrow \infty}\sup\{v(t)\} \leq {M}_1, \,\, \text{and} \quad \lim_{t \rightarrow \infty}\sup\{h(t)\} \leq {M}_2.
\end{eqnarray}
\end{Thm}
\begin{Prf}
From the first equation of system \eqref{NAM} we have
\begin{eqnarray}\label{M}
\begin{array}{l}
\dot v(t) \leq v(t)\left[ a(t)-b(t)v(t) \right].  \nonumber 
\end{array}
\end{eqnarray}
Consider the following equation
\begin{eqnarray}\label{M1}
 \dot U(t) = U(t)\left[ a(t)-b(t)U(t) \right].
\end{eqnarray}
By Lemma (\ref{globally-stable-lemma}), equation \eqref{M1} has a positive $\omega$-periodic solution $U^*(t)$ which is globally asymptotically stable. Let $U(t)$ be any positive solution to \eqref{M1} with $U(0)= v(0)$. Using the comparison method of differential equations (see \cite{Cosner}) we have:
\begin{eqnarray}\label{M11}
v(t) \, \leq \, U(t), \quad t \geq 0
\end{eqnarray}
since $U^*(t)$ is globally attractive, for given $\varepsilon>0$ there exists a $T_0>0$ such that
\begin{equation}\label{M10}
 v(t) < U^*(t) + \varepsilon,  \qquad t>T_0.
\end{equation}
Let $M_1= \max_{t \in \left[T_0, \infty \right)} \{U^*(t)+\varepsilon\}$ and then 
\begin{eqnarray}
 \lim_{t \rightarrow \infty} \sup \{ v(t) \} \leq M_1.
\end{eqnarray}
To prove the second part of the theorem, from the second equation of system \eqref{NAM} we have
\begin{equation}\label{M2}
 \dot h(t) \leq h(t) \left(-R(t) +  \frac{\alpha(t) \left[ v^{*}(t) + \varepsilon - \rho(t) \right]} { \beta(t) + v^{*}(t) + \varepsilon} - \frac{\gamma(t) \beta(t)}{M_1 - \rho(t)} \right)
\end{equation}
where $v^{*}(t)$ is the periodic solution to \eqref{gas}. Integrating the both sides of \eqref{M2} from $T_0$ to $t$ yields:
\begin{equation}\label{M02}
 h(t) \leq h(T_0) \exp \left\{ \int_{T_0}^{t} \left( -R(s) + \frac{\alpha(s) \left[ v^{*}(s) + \varepsilon - \rho(s) \right]} { \beta(s)} - \frac{\gamma(s) \beta(s)}{M_1 - \rho(s)} \right) ds \right\} \nonumber
\end{equation}
and since $v^{*}(t)$ is bounded over $\left[ T_0,t \right]$, then $ \left( -R(t) + \frac{\alpha(t) \left[ v^{*}(t) + \varepsilon - \rho(t) \right]} { \beta(t)} - \frac{\gamma(t) \beta(t)}{M_1 - \rho(t)} \right)$ is also bounded over the same interval. Now we have
\begin{equation}
h(t) \, \leq \, h \left( T_0 \right) \exp\left\{ \eta \right\}
\end{equation}
where
\begin{equation}\label{M201}
\eta \,=\, \int_{T_0}^{t} \left( -R(s) + \frac{\alpha(s) \left[ v^{*}(s) + \varepsilon - \rho(s) \right]} { \beta(s)} - \frac{\gamma(s) \beta(s)}{M_1 - \rho(s)} \right) ds. \nonumber
\end{equation}
Let $M_2 \,=\, \max_{t \in \left[T_0, \infty \right)} \left\{ h \left( T_0 \right) \exp\left( \eta \right) \right\}$.
Then we have
\begin{equation}
  \lim_{t \rightarrow \infty}\sup\{h(t)\} \leq M_2 \nonumber.
\end{equation}
The proof is complete.   
\end{Prf}
The vegetation biomass is always greater than its reserve $v_u$, and then we can conclude that $\lim_{t \rightarrow \infty} \inf\{v(t) \} \geq v_u$. When $v_u = 0$, the following theorem prove that there exists an integer $ m_1 > 0$ such that $\lim_{t \rightarrow \infty} \inf\{v(t) \} \,>\, m_1$.
\begin{Thm}\label{m_1_theorem}
 Consider system \eqref{NAM}. If $A_\omega \left( b \right)\,>\,0$ and $A_\omega \left( a - \frac{c \left( M_1 - \rho \right)}{\beta} M_2 \right)\,>\,0$, where $A_\omega (\cdot)$ defined by \eqref{average}, there exists a positive constant $m_1<M_1$
 such that:
\begin{equation}\label{m1}
 \lim_{t \rightarrow \infty} \inf\{v(t)\} \,>\, m_1. \nonumber
\end{equation}
\end{Thm}
\begin{Prf}
From Theorem \eqref{M_1-M_2-theorem}, there always exists $T_2>0$ such that $0<m_1\leq M_1$. If $v(0)>0$, from the first equation of system \eqref{NAM} we have:
\begin{equation}\label{m10}
 \dot v(t) \geq v(t)\left [a(t)-\frac{c(t) \left( M_1 - \rho(t) \right)}{\beta(t)} M_2 - b(t)v(t) \right]. 
\end{equation}
Consider the following auxiliary equation:
\begin{equation}\label{m11}
 \dot u(t) = u(t)\left [ a(t) - \frac{c(t) \left( M_1 - \rho(t) \right)}{\beta(t)} M_2  - b(t) u(t) \right]
\end{equation}
Let $\bar u(t)$ be the positive $\omega$-periodic solution to \eqref{m11}, $u(t)$ be any positive solution to \eqref{m11} with $u(T_2) = v(T_2)$. From the comparison method of differential equations we have:
\begin{eqnarray}\label{m12}
 v(t) \geq u(t), \qquad t \geq T_2. \nonumber
\end{eqnarray}
From the global attractivity of $\bar u(t)$, there exists $T_3 > T_2$ such that
\begin{eqnarray}
 \frac{\bar u(t)}{2} < v(t), \qquad  t \geq T_3.\nonumber
\end{eqnarray}
Let $m_1=\min\{ \frac{\bar u(t)}{2}\}, \quad  t \geq T_3$, then
\begin{equation}
 \lim_{t \rightarrow \infty} \inf\{v(t) \} > m_1.
\end{equation}
The proof is complete.   
\end{Prf}
\begin{Thm}\label{m_2_theorem}
Consider system \eqref{NAM}. If $A_\omega \left( -R \,+\, \frac{\alpha \left( v^* \,-\, \rho \right)} {\beta \,+\, v^* } - \gamma \left( 1 \,+\, \frac{ \bar \beta } { v^* \,-\, \rho }  \right)  \right) > 0$, where $v^*(t)$ is the unique positive $\omega$-solution to system \eqref{gas}, $A_\omega (\cdot)$ defined by \eqref{average} and $\bar \beta(t) = \beta(t) + \rho(t)$, there exists a positive constant $m_2 \,<\, M_2$ such that
\begin{equation}\label{m2}
 \lim_{t \rightarrow \infty} \inf\{h(t)\} > m_2.
\end{equation}
\end{Thm}
\begin{Prf}
We can choose $\varepsilon_0 \,>\, 0$ such that
\begin{equation}\label{m201}
 A_\omega\left(\Phi_{\varepsilon_0}(t)\right) > 0, \quad \text{and} \quad A_\omega\left(\Phi_\varepsilon(t)\right) > 0
\end{equation}
where $\varepsilon$ is same as in theorem \eqref{M_1-M_2-theorem} and
\begin{equation}\label{m21}
 \Phi_{\varepsilon_0} = -R(t)+ \frac{\alpha(t) \left[ v^*(t)-\varepsilon_0 -\rho(t) \right]} {\beta(t) + v^*(t)-\varepsilon_0 } - \gamma(t) \left[ 1 + \frac{ \bar \beta(t) } { v^*(t) - \varepsilon_0 - \rho(t) }  \right]
\end{equation}
\begin{equation}\label{m22}
\Phi_{\varepsilon} = -R(t)+ \frac{\alpha(t) \left[ v^*(t)-\varepsilon -\rho(t) \right]} {\beta(t) + v^*(t)-\varepsilon } - \gamma(t) \left[ 1 + \frac{ \bar \beta(t) } { v^*(t) - \varepsilon - \rho(t) }  \right].
\end{equation}
From the first equation of system \eqref{NAM}, consider the auxiliary equation with the parameter $\eta \,>\, 0$ :
\begin{equation}\label{m23}
 \dot v(t) = v(t) \left[ a(t)- \left( b(t) + \eta \frac{c(t)}{\beta(t)} \right) v(t) \right].
\end{equation}
If $A_\omega(a) >0, \, \text{and} \, A_{\omega} \left(b(t) + \eta \frac{c(t)}{\beta(t)} \right) > 0$, Lemma \eqref{globally-stable-lemma} tell us that system \eqref{m23} has a unique positive $\omega$-periodic solution $v_\eta (t)$. Let $\bar {v}_\eta (t)$ be any solution to \eqref{m23} with $\bar {v}_\eta (0) = v^*(0)$, then for the above $\varepsilon_0$ and from global asymptotic stability of $v_{\eta}$, there exists $T_5 > T_4$ such that
\begin{equation}
 | \bar {v}_\eta(t) - v_\eta(t)| < \frac{\varepsilon_0}{4}, \quad t > T_5.
\end{equation}
Let $v^*(t)$ be the unique positive $\omega$-periodic solution to \eqref{m23} when $\eta = 0$. From continuity of solution in the parameter $\eta $, there exists $\eta_0 =\eta_0(\varepsilon_0)$ and $t \in  \left[ T_5, T_5 +\omega \right]$ such that
\begin{equation}\label{m4}
 | {v}_\eta(t) - v^*(t)| \, 
 < \, \frac{\varepsilon_0}{2}, \quad t \in  \left[ T_5, T_5 +\omega \right].
\end{equation}
Since $v_\eta (t)$ and $v^*(t)$ are both $\omega$-periodic, then we have
\begin{equation}\label{m40}
 | {v}_\eta(t) - v^*(t)|  < \frac{\varepsilon_0}{2}, \quad t \geq 0, \quad 0 < \eta < \eta_0.
\end{equation}
Let $0 < \eta_1 < \eta_0$. From \eqref{m40} we have
\begin{eqnarray}\label{MM}
 v_{\eta_1} (t) \geq v^*(t)- \frac{\varepsilon_0}{2}, \quad t \geq 0.
\end{eqnarray}
Suppose that \eqref{m2} does not hold, then there exists $Z \in \RR^{2}_+$ such that for all positive solution $(v(t,Z), h(t,Z))$ of system \eqref{NAM} with initial condition $(v(0), h(0))=Z$, we have
\begin{equation}\label{m25}
 h(t,Z) < 2\eta_1, \quad (2\eta_1<\varepsilon_0), \quad t > T_5.
\end{equation}
On other hand theorem \eqref{M_1-M_2-theorem} tell us that, there exists $T_6 > T_5$ such that
\begin{eqnarray}\label{m1_m_2}
 v(t,Z) \leq M_1,  \quad T_6 \,>\, T_5.
\end{eqnarray}
Applying \eqref{m25} and \eqref{m1_m_2} to the first equation of \eqref{NAM}, it follows that for all $t>T_6$
\begin{equation}
 \dot v(t,Z) \geq v(t,Z) \left[ a(t) - \left( b(t) + 2\eta_1 \frac{c(t)}{\beta(t) }\right) v(t,Z) \right].
\end{equation}
Let $u(t)$ be a positive  solution to \eqref{m23} with $\eta=\eta_1$ and $u(T_6) \,=\, v(T_6,Z)$, then we obtain
\begin{equation}
 v(t,Z) \,>\, u(t), \quad t \,>\, T_6
\end{equation}
thus, $u_{\eta_1}(t)$ is periodic solution to \eqref{m23}, so from its global asymptotic stability, for given $\varepsilon = \frac{\varepsilon_0}{2}$, there exists $T_7 > T_6$ such that 
\begin{eqnarray}\label{m267}
 v(t,Z) \,>\, u(t) > v_{\eta_1} (t) - \frac{\varepsilon_0}{2}, \quad  t \geq T_7
\end{eqnarray}
 and hence from \eqref{MM}, and \eqref{m267} we obtain
\begin{equation}\label{m27}
 v(t,Z) > v^*(t) - \varepsilon, \quad t \geq T_7.
\end{equation}
From the second equation of \eqref{NAM} and by applying \eqref{m27} we have
\begin{equation} \label{m29}
 \dot h(t,Z) \geq h(t,Z) \, \min \left \{ \Phi_{\varepsilon_0}(t),  \Phi_{\varepsilon}(t) \right \}
\end{equation}
integrating \eqref{m29} with respect to $t$ from $T_7$ yields
\begin{equation}\label{m290}
 h(t,Z) \geq h(T_7,Z) \exp \left\{ \int_{T_7}^t \min \left \{ \Phi_{\varepsilon_0}(s), \Phi_{\varepsilon}(s) \right \} ds \right \}, \quad t > T_7.
\end{equation}
Thus, from \eqref{m201} it follows that $h(t) \longrightarrow + \infty$, as $t \longrightarrow + \infty$ which is a contradiction to the claim \eqref{m25}. So there must exists $m_2$ such that $\lim_{t \rightarrow \infty} \inf h(t) > m_2$. This complete the proof.
\end{Prf}
\section{Main Results}
\begin{Thm}\label{sufficent-conditions-for-permenancy}\label{sufficent-conditions-for-permenancy-Thm}
If all the assumptions in theorem \eqref{M_1-M_2-theorem}, \eqref{m_1_theorem} and \eqref{m_2_theorem} hold. Then system \eqref{NAM} is permanent.
\end{Thm}
The following theorem gives the necessary condition for permanency of the system \eqref{NAM}.
\begin{Thm}\label{necessary-conditions-for-permenancy}
System \eqref{NAM} is permanent if and only if
\begin{eqnarray}\label{per}
\begin{array}{l}
A_\omega \left ( -R \,+\, \frac{\alpha \left( v^* \,-\, \rho \right) } { \beta \,+\, v^*} - \frac{ \gamma \beta } {v^* \,-\, \rho} \right) > 0
\end{array}
\end{eqnarray}
where $A_\omega(\cdot)$ defined by \eqref{average}, and $v^*(t)$ is the unique positive periodic solution of \eqref{M1}.
\end{Thm}
\begin{Prf}
Suppose that
\begin{equation}\label{main_0}
\begin{array}{l}
 A_\omega \left (-R \,+\, \frac{\alpha \left( v^* \,-\, \rho \right)}{\beta \,+\, v^* } \,-\, \frac{\gamma \beta} {v^* \,-\, \rho} \right) \leq 0
\end{array}
\end{equation}
with this assumption we only need to show that $\lim_{t \longrightarrow \infty} h(t) \,=\, 0$. From the first equation in \eqref{NAM} we have
\begin{equation}
 \dot v(t) \leq v(t)\left[ a(t) \,-\, b(t)v(t) \right].
\end{equation}
Consider the following auxiliary equation
\begin{equation}\label{main_2}
 \dot u(t) = u(t)\left[ a(t)-b(t)u(t) \right].
\end{equation}
Let $v^{*}(t)$ be the unique positive $\omega$-periodic solution to \eqref{main_2} and $u(t)$ be any positive solution to \eqref{main_2}. From global asymptotic stability of $v^{*}(t)$ and using the comparison method, for any positive $\varepsilon < 1$, there exists $T^{(1)}>0$ such that
\begin{eqnarray}\label{per202}
u(t) < v^*(t) +\varepsilon, \quad t \geq T^{(1)}.
\end{eqnarray}
By comparison method and from \eqref{per202} we have
\begin{eqnarray}\label{per2}
 v(t) \leq v^*(t)+\varepsilon, \quad t \geq T^{(1)}.
\end{eqnarray}
It follows that
\begin{eqnarray}\label{per01}
 \frac{1}{\omega}\int_{0}^\omega \left (-R(t) + \alpha(t)\, \frac{ v^{*}(t) + \varepsilon -\rho(t)}{\beta(t) + v^{*}(t) + \varepsilon} - \frac{\gamma(t) \left( v^{*}(t) + \varepsilon - \beta(t) \right) } {v^{*}(t) + \varepsilon -\rho(t)} 
\right)dt \, < \, -\varepsilon_0.
\end{eqnarray}
From \eqref{per2}  we also have
\begin{eqnarray}\label{per1}
 \frac{1}{\omega}\int_{0}^\omega \left (-R(t) + \alpha(t)\, \frac{ v(t) -\rho(t)}{\beta(t) + v(t)} - \frac{\gamma(t) \left[ v(t) - \beta(t) \right] } {v(t) -\rho(t)} 
\right)dt \, < \, -\varepsilon_0.
\end{eqnarray}
We will show that for the above $\varepsilon$, there exists a $T^{(2)}>T^{(1)}$ such that $h(T^{(2)})<\varepsilon$, otherwise from \eqref{per1} we have
\begin{align}
  \varepsilon \leq h(t)  
\leq  h(T^{(1)})\exp\left[ \int_{T^{(1)}}^t \left (-R(s)+ \alpha(s)\, \frac{  v(s) -\rho(s) }{\beta(s) + v(s)} 
\right)ds \right] \longrightarrow 0 
\end{align}
as $ \, t \longrightarrow +\infty$, which is a contradiction.\\
Let
\begin{eqnarray}
  M(\varepsilon)= \max_{t \in [0,\omega]} \left\{ -R(t)+ \frac{\alpha(t)\left(v^*(t) + \varepsilon -\rho(t) \right)} {\beta(t)+\left(v^*(t)+ \varepsilon \right)} - \frac{\gamma(t) \left[ v^*(t) + \varepsilon - \beta(t) \right] } { v^*(t) + \varepsilon - \rho(t)} \right\}
\end{eqnarray}
which is bounded for $\varepsilon \in [0,1]$. Now we only need to show
\begin{eqnarray}\label{main}
h(T^{(3)}) \leq \varepsilon \exp\left\{ M(\varepsilon)\omega\right\}
\end{eqnarray}
otherwise there exists a $T^{(3)}>T^{(2)}$ such that
\begin{eqnarray}
h(T^{(3)})>\varepsilon \exp\left\{ M(\varepsilon)\omega\right\}.
\end{eqnarray}
By the continuity of $h(t)$ there exists $T^{(4)} \in \left(T^{(2)},T^{(3)}\right)$ such that $h(T^{(4)})=\varepsilon$ and $h(t)>\varepsilon$ for $t \in \left(T^{(4)},T^{(3)} \right)$.\\
Let $P$ be the nonnegative integer such that $T^{(3)} \in \left(T^{(4)} + P\omega,T^{(4)}+(P+1)\omega \right)$, from \eqref{per1} we have
\begin{eqnarray}
\begin{array}{l}
\varepsilon \exp\{ M(\varepsilon)\omega\} \quad < \quad h(T^{(3)})
 \\
 \qquad \leq \quad h(T^{(4)}) \exp\left\{ \int_{T^{(4)}}^{T^{(3)}} \left (-R(t)+ \frac{\alpha(t) \left( v(t) -\rho (t) \right)} {\beta(t) + v(t)} - \frac{\gamma(t) \left[ v(t) - \beta(t) \right] } { v(t) -\rho(t)} \right)dt \right\} 
\\
 \qquad = \quad \varepsilon \exp\left\{ \int_{T^{(4)}}^{T^{(4)}+P\omega} + \int_{T^{(4)}+P\omega}^{T^{(3)}}  \left (-R(t)+ \frac{\alpha(t) \left( v(t) -\rho (t) \right)} {\beta(t) + v(t)} - \frac{\gamma(t) \left[ v(t) - \beta(t) \right] } { v(t) -\rho(t)} \right)dt \right\} \nonumber
\\
 \qquad \leq \quad \varepsilon \exp\left\{\int_{T^{(4)} + P\omega}^{T^{(3)}}  \left (-R(t)+ \frac{\alpha(t) \left( v(t) -\rho (t) \right)} {\beta(t) + v(t)} - \frac{\gamma(t) \left[ v(t) - \beta(t) \right] } { v(t) -\rho(t)} \right)dt \right\} 
\\
\qquad < \quad \varepsilon \exp\{ M(\varepsilon)\omega\}
\end{array}
\end{eqnarray}
which is contradiction. It follows that $h(t) \longrightarrow 0 $ as $t \longrightarrow +\infty$, which contradict that $h(t)$ is permanent, so \eqref{main_0} must be false i.e. \eqref{per} holds. This complete the proof.   
\end{Prf}
The following theorem gives the sufficient conditions for \eqref{NAM} having a globally asymptotically stable solutions.   Let
\begin{eqnarray}\label{Gamma_epsilon}
\Gamma = \left\{ \left(v,h \right) \in \RR^2 |\, m_1\, \leq\, v \,\leq\, M_1; \, m_2 \,\leq\, h \,\leq\, M_2 \right\}.
\end{eqnarray}
Suppose that
\begin{eqnarray} \label{Thm-Conditions}
\left.
\begin{array}{l}
\, \inf_{t \in \RR} \left\{ \,   b(t) \,+\, c(t) \,\frac{  \rho(t) \left( m_1^{} \right) - \left( m_1^{} \hat v \,-\, \rho(t) \beta(t) \right) \hat h } { M_1^{} \hat v(t) \left(\beta(t) \,+\,  M_1^{} \right) (\beta(t) \,+\, \hat v(t)) } \, - \frac{\alpha (t) \beta (t) \,+\, \alpha (t) \rho (t)}{\left( \beta (t) \,+\, m_1^{} \right)  \left( \beta (t) \,+\, \hat v(t) \right)}\, \right.
\\
\qquad \qquad \left. - \frac{\gamma (t) \beta (t) \,-\, \gamma(t) \rho(t)}{ \left( m_1^{} \,-\, \rho(t) \right)\left( \hat v(t) \,-\, \rho(t) \right) }   \right\} \, > \, 0 ,
\\
\text{and} \\
\,  \inf_{t \in \RR} \left\{ \, c(t) \,\frac{\hat v(t) \left[ m_1^{} \left( \beta(t) \,+\, \hat v(t) \right) \, -\, \hat v(t) \rho (t) \left( \beta(t) \,+\, \rho (t) \hat v(t) \right) \right]}{M_1^{} \hat v(t) \left( \beta(t) \,+\, M_1^{} \right) M_1^{} \hat v(t) \left( \beta(t) \,+\, \hat v(t) \right) } \right\} \, > \, 0 ,
\end{array} 
\right.
\end{eqnarray}
where $m_1^{}\,\, \text{and} \,\,M_1^{}$ defined by \eqref{Gamma_epsilon}.
\begin{Thm}
 Let $(\hat v(t), \hat h(t))$ be a bounded positive solution of \eqref{NAM}. If \eqref{Thm-Conditions} holds, and if
\begin{eqnarray}
\left.
\begin{array}{l}
A_\omega \left( a - \frac{c \left( M_1 -\rho \right)}{\beta} M_2 \right) \, > \, 0, \,\, 
 A_\omega \left (-R + \frac{\alpha \left( \hat v \,-\, \rho \right)}{ \beta \,+\, \hat v } - \frac{\gamma \beta} {\hat v \,-\, \rho } \right)\, > \, 0 ,
\end{array} 
\right.
\end{eqnarray}
where $A_\omega (\cdot)$ defined by \eqref{average}. Then $(\hat v(t), \hat h(t))$ is globally asymptotically stable.
\end{Thm}
\begin{Prf}
Let $(v(t), h(t))$ be any solution of \eqref{NAM} with $(v(t_0) > 0, h(t_0) > 0)$. The set $\Gamma$ defined by \eqref{Gamma_epsilon} is an ultimately bounded region of \eqref{NAM}, so there exists a $T^{(5)} > 0$ such that $(v(t), h(t)) \in \Gamma$ and $(\hat v(t), \hat h(t)) \in \Gamma$ for all $t \geq t_0 + T^{(5)}$.
Consider a Lyapunov function defined by
\begin{equation}\label{lyapn}
  X(t) = | \ln v(t)- \ln \hat v(t)| + | \ln h(t)- \ln \hat h(t)|;\quad t > t_0.
\end{equation}
Calculating the right upper derivative of \eqref{lyapn} produces
\begin{eqnarray}
\begin{array}{l}
 X'(t)  = \text{sgn} \left\{ \ln v(t) - \ln \hat v (t) \right\} \left[\frac{v'(t)}{v(t)} - \frac{(\hat v (t))'}{\hat v (t)} \right] + \text{sgn} \{\ln h(t) - \ln \hat h (t)\}\left[\frac{h'(t)}{h(t)} - \frac{(\hat h (t))'}{\hat h (t)} \right]
\\
\\
\, = \, \text{sgn} \left\{ \ln v(t) - \ln \hat v (t) \right\} \left[ -b(t) (v(t) - \hat v(t)) - c(t) \left( \frac{v(t) \,-\, \rho (t) ) h(t)}{v(t) (\beta(t) \,+\, v(t))} - \frac{\hat v(t) \,-\, \rho (t) ) \hat h(t)}{\hat v(t) (\beta(t) \,+\, \hat v(t))} \right) \right]
\\
\\
\, + \, \text{sgn} \left\{ \ln h(t) - \ln \hat h (t)\right\} \left[ \frac{\alpha (t) (v(t) \,-\, \rho (t))}{\beta(t) \,+\, v(t)} - \frac{\gamma(t) (\beta(t) \,+\, v(t))} { v(t) \,-\, \rho (t) } - \frac{\alpha (t) ( \hat v(t) \,-\, \rho (t))}{\beta(t) \,+\, \hat v(t)} + \frac{\gamma(t) \left(\beta(t) \,+\, \hat v(t) \right)} { \hat v(t) - \rho (t) } \right] \nonumber
\end{array}
\end{eqnarray}
\begin{eqnarray}
\begin{array}{l}
\, \leq \, \text{sgn} \left\{ \ln v(t) - \ln \hat v (t) \right\} \left\{ \left( -b(t) \,-\, c(t) \frac{ \rho(t) \left(v(t) \,+\, \hat v(t)\right) \,-\, \left( v(t) \hat v(t) \,-\, \rho (t) \beta(t) \right) \hat h(t) }{ v(t) \hat v(t) \left( \beta(t) \,+\, v(t)\right) \left( \beta(t) \,+\, \hat v(t) \right) } \right) [v(t) \,-\, \hat v(t)]  \right.
\\
\\
\qquad \left. -\, c(t) \left(  \frac{\hat v(t) \left( v[\beta(t) \,+\, \hat v(t)] \,-\, \hat v(t)\rho(t) [ \beta(t) \,+\, \rho(t) \hat v(t)] \right)}{v(t) \hat v(t) \left( \beta(t) \,+\, v(t)\right) \left( \beta(t) \,+\, \hat v(t) \right) }\right) \left[ h(t) \,-\, \hat h(t) \right] \right\}
\\
\\
\qquad + \,\, \text{sgn} \left\{ \ln h(t) - \ln \hat h (t) \right\} \left\{ \frac{\alpha(t) \beta(t) \,+\, \alpha (t) \rho(t)}{\left(\beta (t) \,+\, v(t) \right) \left(\beta (t) \,+\, \hat v(t) \right) } \,+\, \frac{\gamma(t) \beta(t) \,-\, \gamma(t) \rho(t)}{ \left( v(t) \,-\, \rho(t) \right) \left( \hat v(t) \,-\, \rho(t) \right) } \right\} [v(t) \,-\, \hat v(t)] \nonumber
\end{array}
\end{eqnarray}
\begin{eqnarray}
\begin{array}{l}
\, \leq \, \left[ -\, b(t) \,-\, c(t) \frac{ \rho(t) \left( m_1^{} \,+\, \hat v(t)\right) \,-\, \left( m_1^{} \hat v(t) \,-\, \rho (t) \beta(t) \right) \hat h(t) }{ M_1^{} \hat v(t) \left( \beta(t) \,+\, m_1^{} \right) \left( \beta(t) \,+\, \hat v(t) \right) } \,+ \, \frac{\alpha(t) \beta(t) \,+\, \alpha (t) \rho(t)}{\left(\beta (t) \,+\, m_1^{} \right) \left(\beta (t) \,+\, \hat v(t) \right) } \right.
\\
\\
\quad \left. \,+\, \frac{\gamma(t) \beta(t) \,-\, \gamma(t) \rho(t)}{ \left( m_1^{} \,-\, \rho(t) \right) \left( \hat v(t) \,-\, \rho(t) \right) } \right] [v(t) \,-\, \hat v(t)] \nonumber
\,-\, c(t) \left[  \frac{\hat v(t) \left( m_1^{}[\beta(t) \,+\, \hat v(t)] \,-\, \hat v(t)\rho(t) [ \beta(t) \,+\, \rho(t) \hat v(t)] \right)}{M_1^{} \hat v(t) \left( \beta(t) \,+\, M_1^{}\right) \left( \beta(t) \,+\, \hat v(t) \right) }\right] \left| h(t) \,-\, \hat h(t) \right|
\\
\\
= - \left[ b(t) + c(t) \frac{ \rho(t) \left( m_1^{} \,+\, \hat v(t)\right) \,-\, \left( m_1^{} \hat v(t) \,-\, \rho (t) \beta(t) \right) \hat h(t) }{ M_1^{} \hat v(t) \left( \beta(t) \,+\, m_1^{} \right) \left( \beta(t) \,+\, \hat v(t) \right) } - \frac{\alpha(t) \beta(t) \,+\, \alpha (t) \rho(t)}{\left(\beta (t) \,+\, m_1^{} \right) \left(\beta (t) \,+\, \hat v(t) \right) } \right. \nonumber
\\
\\
\quad \left. - \frac{\gamma(t) \beta(t) \,-\, \gamma(t) \rho(t)}{ \left( m_1^{} \,-\, \rho(t) \right) \left( \hat v(t) \,-\, \rho(t) \right) } \right] |v(t) - \hat v(t)| \,-\, c(t) \left[  \frac{\hat v(t) \left( m_1^{}[\beta(t) \,+\, \hat v(t)] \,-\, \hat v(t)\rho(t) [ \beta(t) \,+\, \rho(t) \hat v(t)] \right)}{M_1^{} \hat v(t) \left( \beta(t) \,+\, M_1^{}\right) \left( \beta(t) \,+\, \hat v(t) \right) }\right] \left| h(t) \,-\, \hat h(t) \right|
\end{array}
\end{eqnarray}
and since
\begin{eqnarray}
 \begin{array}{l}
  \inf_{t \, \in\, \RR} \left\{ b(t) \,+\, c(t) \frac{ \rho(t) \left( m_1^{} \,+\, \hat v(t)\right) \,-\, \left( m_1^{} \hat v(t) \,-\, \rho (t) \beta(t) \right) \hat h(t) }{ M_1^{} \hat v(t) \left( \beta(t) \,+\, m_1^{} \right) \left( \beta(t) \,+\, \hat v(t) \right) } \,-\, \frac{\alpha(t) \beta(t) \,+\, \alpha (t) \rho(t)}{\left(\beta (t) \,+\, m_1^{} \right) \left(\beta (t) \,+\, \hat v(t) \right) } \right.
 \\
 \qquad -\, \left. \frac{\gamma(t) \beta(t) \,-\, \gamma(t) \rho(t)}{ \left( m_1^{} \,-\, \rho(t) \right) \left( \hat v(t) \,-\, \rho(t) \right) } \right\} \,>\, 0 \nonumber
 \\
 \text{and} \quad 
 \inf_{t\,\in\,\RR}  \left\{ c(t)  \frac{\hat v(t) \left( m_1^{}[\beta(t) \,+\, \hat v(t)] \,-\, \hat v(t)\rho(t) [ \beta(t) \,+\, \rho(t) \hat v(t)] \right)}{M_1^{} \hat v(t) \left( \beta(t) \,+\, M_1^{}\right) \left( \beta(t) \,+\, \hat v(t) \right) }\right\} \, > \, 0 ,
 \end{array}
\end{eqnarray}
there exists a positive constant $\mu$ such that
\begin{eqnarray}\label{lyapn1}
 X'(t) \leq -\mu \left(\left|v(t)-\hat v(t) \right|+ \left|h(t)-\hat h(t)\right| \right), \qquad    t\,>\,t_0 \,+\, T^{(5)}.
\end{eqnarray}
Integrating both sides of \eqref{lyapn1} with respect to $t$ from $T+t_0$ yields
\begin{eqnarray}
 X(t)-X(t_0 \,+\, T^{(5)}) \leq -\mu \int_{t_0 \,+\, T^{(5)}}^t\left( |v(s)-\hat v(s)|+ \left|h(s)-\hat h(s)\right| \right)ds \nonumber
\\
X(t)+\mu \int_{t_0 \,+\, T^{(5)}}^t\left(|v(s)-\hat v(s)| + \left|h(s)-\hat h(s) \right| \right)ds \leq X(t_0 \,+\, T^{(5)}) \nonumber
\end{eqnarray}
which yields
\begin{eqnarray}
 \mu \int_{t_0 \,+\, T^{(5)}}^t \left(|v(s)-\hat v(s)| + \left|h(s)-\hat h(s) \right| \right)ds \leq X(t_0 \,+\, T^{(5)}) \nonumber
\end{eqnarray}
from here we can conclude that
\begin{eqnarray}
 \int_{t_0 \,+\, T^{(5)}}^t \left(|v(s)-\hat v(s)| + \left|h(s)-\hat h(s) \right| \right)ds \leq \mu^{-1}X(t_0 \,+\, T^{(5)}) < + \infty \nonumber
\end{eqnarray}
and hence $|v(t)-\hat v(t)|+ \left|h(t)-\hat h(t)\right| \in L^{1}\left([t_0 \,+\, T^{(5)}, +\infty)\right)$. Since $\hat v (t), \hat h (t)$ is bounded and $v(t), h(t)$ is ultimately bounded, this implies that they all have bounded derivatives for $t\geq t_0+T$. It follows that $|v(t)-\hat v(t)|+ \left| h(t)-\hat h(t)\right|$ is uniformly continuous on $[t_0 \,+\, T^{(5)}, +\infty)$ and hence using By Barb\u{a}lat's Lemma \cite{Barb} one can conclude that
\begin{eqnarray}
 \lim_{t\leftarrow + \infty} \left( \left| v(t)-\hat v(t) \right| + \left| h(t)-\hat h(t) \right| \right) \, = \,0.
\end{eqnarray}
This shows that $\left( \hat v(t), \hat h(t) \right)$ is a globally asymptotically stable solution to \eqref{NAM}. The proof is complete.
\end{Prf}
Let $ \left( v(\omega, t_0,(v_0,h_0)) \, ,\\ h(\omega, t_0,(v_0,h_0)) \right) $ be the solution of \eqref{NAM} through the point $\left( t_0,(v_0,h_0) \right)$.
\begin{Thm}
Consider the set \eqref{Gamma_epsilon}. If $\, \, A_\omega \left( -R + \frac{\alpha \left( \bar v + \varepsilon - \rho \right)} {\beta} - \frac{\gamma \beta} {M_1 - \rho} \right)\, > \, 0, \,\, \left. A_\omega \left( \frac{\gamma} {\alpha \left( M_1  - \rho \right)} \right) \, > \, 0 \right.\\ \text{and} \quad A_\omega \left( a - \frac{cM_2}{\beta} \right) > 0$ , where $A_\omega(\cdot)$ is defined by \eqref{average}, and $\bar v$ is the periodic solution of \eqref{gas}, then system \eqref{NAM} has at least one positive $\omega$-periodic solution $(v^*,h^*)$ which lies in $\Gamma$.
\end{Thm}
\begin{Prf}
Define a Poincaré$\acute{e}$ shift operator $\sigma: \RR_{+}^2 \longrightarrow \RR_{+}^2$, by
\begin{equation}
 \sigma \left( (v_0,h_0) \right) = \left( v(\omega, t_0,(v_0,h_0)),\, h(\omega, t_0,(v_0,h_0)) \right), \qquad (v_0,h_0)\, \in \, \RR_{+}^2
\end{equation}
Thus, from theorem \eqref{sufficent-conditions-for-permenancy} we know that $\Gamma$ is positively invariant of \eqref{NAM}, and hence $\sigma(\Gamma)=\Gamma$. Since the solution of \eqref{NAM} is continuous with respect to the initial value $(t_0,(v_0,h_0))$, then the operator $\sigma$ is continuous. From the definition of the set $\Gamma$ we know that it is bounded, closed, and convex. By Brouwer's fixed-point Theorem \cite{Gerard1959} we know that $\sigma$ has at least one fixed point
\begin{eqnarray}
 (v^*,h^*) \,=\, \sigma((v^*,h^*)) \,=\, \left( v(\omega,t_0,(v^*,h^*),\, h(\omega,t_0,(v^*,h^*)\right)
\end{eqnarray}
So there exist at least one periodic solution $(v^*,h^*)\in \Gamma$. The proof is complete.
\end{Prf}

\bibliographystyle{unsrtnat}
\bibliography{myref1} 





\end{document}